\documentclass[a4paper]{article}

\usepackage{amsmath}
\usepackage{amsfonts}
\usepackage{amssymb}
\usepackage{amsthm}

\title{Symmetric coinvariant algebras and local Weyl modules at a double point}
\author{Toshiro Kuwabara \\ Department of Mathematics, Graduate School
of Science, Kyoto University, Kyoto 606-8502, Japan \\
toshiro@math.kyoto-u.ac.jp}

\newtheorem{definition}{Definition}
\newtheorem{proposition}[definition]{Proposition}
\newtheorem{theorem}[definition]{Theorem}
\newtheorem{corollary}[definition]{Corollary}
\newtheorem{lemma}[definition]{Lemma}

\newcommand{\refprop}[1]{Proposition~\ref{#1}}
\newcommand{\refthm}[1]{Theorem~\ref{#1}}

\newcommand{\reflemma}[1]{Lemma~\ref{#1}}
\newcommand{\refeq}[1]{(\ref{#1})}
\newcommand{\refsec}[1]{Section~\ref{#1}}

\newcommand{\C}{{\mathbb C}}
\newcommand{\N}{{\mathbb N}}
\newcommand{\Z}{{\mathbb Z}}

\newcommand{\Slr}{\mathfrak{sl}_{r+1}}
\newcommand{\g}{\mathfrak{g}}
\newcommand{\Ind}{\operatorname{Ind}}
\newcommand{\Hom}{\operatorname{Hom}}
\newcommand{\Dim}{\mathrm{dim}}
\newcommand{\Ker}{\mathrm{Ker}}
\renewcommand{\Im}{\mathrm{Im}}
\newcommand{\gr}{\operatorname{gr}}
\newcommand{\Span}{\mathrm{span}}

\newcommand{\A}{\mathcal{A}}

\newcommand{\M}{\mathcal{M}}
\newcommand{\X}{\mathcal{X}}

\begin{document}
\begin{center}
{\large\bf Symmetric coinvariant algebras and local Weyl modules at a
 double point} 

\vspace*{0.7cm}
Toshiro Kuwabara 

\vspace*{0.15cm}
Department of Mathematics, Graduate School
of Science, Kyoto University, Kyoto 606-8502, Japan \\
Email: \verb+toshiro@math.kyoto-u.ac.jp+
\end{center}
\begin{abstract}
 The symmetric coinvariant algebra 
 $\C[x_1, \dots, x_n]_{S_n}$
 is the quotient algebra of the polynomial ring by
 the ideal generated by symmetric polynomials vanishing at the
 origin. It is known that the algebra is isomorphic to
 the regular representation of $S_n$. 

 Replacing $\C[x]$ with $A = \C[x,y]/(xy)$, we introduce another
 symmetric coinvariant algebra 
 $A^{\otimes n}_{S_n}$ and 
 determine its $S_n$-module structure. As an application, 
 we determine
 the $\Slr$-module structure of the local Weyl module at a double point for
 $\Slr \otimes A$.
\end{abstract}


\section{Introduction}
\label{sec:introduction}
The symmetric group $S_n$ acts on the polynomial ring of $n$ variables
$\C[x_1, \dots, x_n]$. Let $\C[x_1, \dots, x_n]^{S_n}_{+}$ be the set
of symmetric polynomials vanishing at the origin
$x_1= \dots =x_n=0$. For an algebra $R$ and a subset $S$ of $R$,
let $\langle S \rangle_{R}$ be the ideal of $R$ generated by
$S$. The classical 
symmetric coinvariant algebra $\C[x_1, \dots , x_n]_{S_n}$ is defined
as the quotient algebra:
\[
 \C[x_1, \dots, x_n]_{S_n} = \C[x_1, \dots, x_n] / 
 \langle \C[x_1, \dots, x_n]^{S_n}_{+} \rangle_{\C[x_1, \dots, x_n]}.
\]
It is known that this is isomorphic to the regular
representation of $S_n$ as an $S_n$-module (\cite{Chev}).

The symmetric group $S_n$ acts diagonally on the polynomial ring of
$2n$ variables $\C[x_1, \dots, x_n, y_1, \dots, y_n]$, i.e.
\[
 \sigma P(x_1, \dots, x_n, y_1, \dots, y_n)
 = P(x_{\sigma(1)}, \dots, x_{\sigma(n)}, y_{\sigma(1)}, \dots,
 y_{\sigma(n)})
\]
for $\sigma \in S_n$ and $P \in \C[x_1, \dots, x_n, y_1, \dots, y_n]$.
Recently, Haiman defined the diagonal symmetric coinvariant algebra
\begin{multline*}
 \C[x_1, \dots, x_n, y_1, \dots, y_n]_{S_n} = \\
 \C[x_1, \dots, x_n, y_1, \dots, y_n]/
 \langle \C[x_1, \dots, x_n, y_1, \dots, y_n]^{S_n}_{+} 
 \rangle_{\C[x_1, \dots, x_n, y_1, \dots, y_n]}
\end{multline*}
where $\C[x_1, \dots, x_n, y_1, \dots, y_n]^{S_n}_{+}$ is
symmetric polynomials vanishing at the origin. He determined its 
$S_n$-module structure in the form
\[
 \C[x_1, \dots, x_n, y_1, \dots, y_n]_{S_n} \simeq
 \C PF_n \otimes L_{(1^n)}
\]
where $PF_n$ is the set of parking functions, functions
from $\{1, \dots, n\}$ to
itself satisfying some condition, 
$\C PF_n$ is the vector space spanned by $PF_n$, and 
$L_{(1^n)}$ is the sign representation of $S_n$ (\cite{Haiman}).
In generally, for a partition $\lambda$ we denote by $L_{\lambda}$
the irreducible representation of $S_n$ corresponding to $\lambda$.

Let $M$ be an affine variety over $\C$ and let $A$ be its coordinate ring.
The symmetric coinvariant algebra $A^{\otimes n}_{S_n}$ 
is introduced by Feigin and Loktev in \cite{FL}.
The symmetric group 
$S_n$ acts on $A^{\otimes n}$, the $n$-th tensor product of $A$. 
Fix a base point $0$ on $M$. Let $(A^{\otimes n})^{S_n}_{+}$ be the set
of symmetric elements vanishing at the point $(0, \dots, 0)$.
The symmetric coinvariant algebra is defined as
\[
 A^{\otimes n}_{S_n} = A^{\otimes n} / \langle (A^{\otimes n})^{S_n}_{+} \rangle_{A^{\otimes n}}.
\]

This representation is used for the study of the
structure of $(\Slr \otimes A)$-module $W_M(\{0\}_{\lambda})$ called
the local Weyl module.
Let $\Slr = \mathfrak{n}_{+} \oplus \mathfrak{h} \oplus
\mathfrak{n}_{-}$ be the triangular decomposition of $\Slr$, and
let $\lambda \in \mathfrak{h}^{*}$ be a dominant integrable weight.
The local Weyl module $W_M(\{0\}_{\lambda})$ is the maximal 
$\Slr$-integrable 
$(\Slr \otimes A)$-module generated by a cyclic vector $v_0$
with the following properties:
\[
 (\mathfrak{n}_{+} \otimes P) v_0 = 0, 
 \quad (h \otimes P) v_0 = \lambda(h) P(0) v_0
 \quad \mbox{for all $P \in A$, $h \in \mathfrak{h}$}.
\]
This definition was first
given by Chari and Pressley in \cite{ChariPressley} for $A=\C[x]$
and then generalized by Feigin and Loktev in \cite{FL}.

Let $V_{r+1} = \C^{r+1}$ be the vector representation of $\Slr$ and let
$\omega_1$ be the highest weight of $V_{r+1}$. In \cite{FL},
Feigin and Loktev show that there is an isomorphism of $\Slr$-modules
$:$
\[
 W_M(\{0\}_{n\omega_1}) \simeq 
 \left(
 V_{r+1}^{\otimes n} \otimes A^{\otimes n}_{S_n}
 \right)^{S_n}.
\]
This isomorphism gives us the connection between the $S_n$-module
structure of the symmetric coinvariant algebra and 
the $\Slr$-module structure of the local Weyl module.

In this paper, We consider the case of  $A=\C[x,y]/(xy)$. 
In this case, the corresponding affine variety $M$ has the double point
$0$. We consider the symmetric coinvariant algebra and the
local Weyl module at the double point $0$. Our main result is
\begin{theorem}
 \label{thm:3}
 We have the following isomorphism of $S_n$-modules$:$
 \[
  A^{\otimes n}_{S_n} \simeq \C[S_n] \oplus (n-1) \Ind^{S_n}_{S_2} L_{(1,1)}.
 \]
 where $L_{(1,1)}$ is the sign representation of $S_2$.
\end{theorem}

As a corollary of \refthm{thm:3}, we determine the structure of 
the local Weyl module $W_M(\{0\}_{n\omega_1})$.
\begin{proposition}
 For $n \in \Z_{\geq 0}$, we have
 \[
  W_M(\{0\}_{n \omega_1}) \simeq V_{r+1}^{\otimes n} \oplus
 (n-1) \left(
 V_{r+1}^{\otimes n-2} \otimes \wedge^{2} V_{r+1}
 \right)
 \] as an $\Slr$-module.
\end{proposition}

Let us give a sketch of the proof of \refthm{thm:3}.
We introduce a generalization of the symmetric coinvariant algebra
$R^n_{i,j}$. Let 
$e_1$, $\dots$, $e_{n}$ be the elementary symmetric polynomials
of variables $x_1$, $\dots$, $x_n$, and $f_1$, $\dots$, $f_{n}$ 
these of $y_1$, $\dots$, $y_n$. 
For $I = \{k_1, \dots, k_i\} \subset \{1, \dots, n\}$,
let $x_I = x_{k_1} \dots x_{k_i}$
and let $y_I = y_{k_1}\dots y_{k_n}$. The algebra $R^n_{i,j}$ is
defined as
\[
 R^n_{i,j} = A^{\otimes n} /
 \langle
 e_1, \dots, e_{i-1}, x_I\:(|I|=i), f_1, \dots, f_{j-1}, y_J\:(|J|=j)
 \rangle_{A^{\otimes n}}.
\]
Clearly we have $R^n_{n,n} = A^{\otimes n}_{S_n}$. By using the same method as
that in \cite{GP}, we can determine the $S_n$-module structure of
$R^n_{i,j}$ for $i,j \geq 1$, $i+j \leq n+1$:
\begin{equation}
 \label{eq:2}
  R^n_{i,j} \simeq \Ind^{S_n}_{S_{n-i-j+2}} L_{(n-i-j+2)}.
\end{equation}

Next, we introduce the decreasing filtration $\{F^p A^{\otimes n}\}_{0 \leq p \leq n}$ of
$A^{\otimes n}$ given by
\[
 F^p A^{\otimes n} = \sum_{k_1 < \dots < k_p} y_{k_1} \dots y_{k_p} A^{\otimes n}.
\]
Let $\{F^p A^{\otimes n}_{S_n}\}_{0 \leq p \leq n}$ be its induced
filtration on $A^{\otimes n}_{S_n}$. 
For $1 \leq i \leq n-1$, We have the following exact sequence 
\begin{equation}
 \label{eq:3}
  0 \rightarrow \gr^i A^{\otimes n}_{S_n} \rightarrow
 R^n_{n-i,i+1} \rightarrow R^n_{n-i,i} \rightarrow 0,
\end{equation}
where $\gr A^{\otimes n}_{S_n}$ is the graded module associated to 
$\{F^p A^{\otimes n}_{S_n}\}_{0 \leq p \leq n}$.
By combining \refeq{eq:2} and \refeq{eq:3}, we obtain \refthm{thm:3}.

The paper is organized as follows. In \refsec{sec:preliminaries}, we recall
basic definitions and notations. In \refsec{sec:symm-coinv-algebra}, 
the symmetric coinvariant algebra is defined. In
\refsec{sec:some-gener-symm}, we introduce the generalization of the
symmetric coinvariant algebra and prove \refeq{eq:2}. 
In \refsec{sec:structure-r_n}, we prove \refeq{eq:3} and then \refthm{thm:3}. 
In \refsec{sec:local-weyl-module}, we review the definition of the
local Weyl module and determine its structure.

\subsection*{Acknowledgments.}
This research is partially supported by Grant-in-Aid for JSPS Research
Fellows No.16-1089.
The author is deeply grateful to Masaki Kashiwara,
Yoshihiro Takeyama and Boris Feigin for 
useful discussions, Sergai Loktev for comments about 
\refsec{sec:local-weyl-module}. 
The author also thanks his advisor Tetsuji
Miwa for reading the manuscript and for his kind encouragement.

\section{Preliminaries}
\label{sec:preliminaries}
In this section we review some definitions and notations in the
representation theory of symmetric groups and symmetric polynomials.

Let $S_n$ be the $n$-th symmetric group. For each partition $\lambda$ of $n$, 
let $L_{\lambda}$ be the irreducible representation of $S_n$
corresponding to $\lambda$.

For a finite group $G$, we denote by $\C[G]$ its group ring. For a $G$-module
$L$ and a subgroup $H$ of $G$, $L^H$ is the subspace of $H$-invariants. 
For an $H$-module $L$, we denote the induced module of $L$ by
$\Ind^{G}_{H} L = \C[G] \otimes_{\C[H]} L$.

The following lemma easily follows from the semi-simplicity of 
representations of $S_n$.
\begin{lemma}
 \label{lemma:6}
 If an $S_n$-module $L$ has a filtration invariant
 under the action of $S_n$, we have an isomorphism 
 $L \simeq \gr L$ of $S_n$-modules where $\gr L$ is the graded
 module associated with the filtration of $L$.
\end{lemma}

Let $V$ be a vector space. We denote by $V^{\otimes n}$ the $n$-th tensor
product of $V$, and by $S^n(V)$ the $n$-th symmetric tensor product. 
For a subset $B \subset V$, we denote by $\Span_{\C} B$ the subspace 
spanned by $B$.

For a set of indeterminates $S$, $e_k(S)$ is the $k$-th elementary
symmetric polynomial of variables $S$. We write $e_1$, $\dots$, $e_n$
for $e_i = e_i(\{x_1, \dots, x_n\})$, and $f_1$, $\dots$, $f_n$ for
$f_i = e_i(\{y_1, \dots, y_n\})$. When the number of variables matters,
we use notations $e_1^{(n)}$, $\dots$, $e_n^{(n)}$ and $f_1^{(n)}$, $\dots$,
$f_n^{(n)}$.

For a set of indices $I = \{j_1, \dots, j_i\} \subset \{1, \dots, n\}$
with $j_1 < \dots < j_i$,
we define
$x_I = x_{j_1} \dots x_{j_i}$.
We also define $y_I$ similarly. 

For an algebra $R$ and a subset $S$ of $R$, 
let $\langle S \rangle_{R}$ be the ideal of $R$
generated by $S$.

\section{The symmetric coinvariant algebra $R_n$}
\label{sec:symm-coinv-algebra}
First we review the classical symmetric coinvariant
algebra $\C[x_1, \dots, x_n]_{S_n}$. The symmetric group $S_n$ acts on
$\C[x_1, \dots, x_n]$. Therefore we can think of the ring of 
symmetric polynomials $\C[x_1, \dots, x_n]^{S_n}$, and we set
\[
 \C[x_1, \dots, x_n]^{S_n}_{+} = 
 \{ P \in \C[x_1, \dots, x_n]^{S_n} \; \vert \;
 P(0, \dots, 0) = 0\}.
\]
Consider the quotient algebra
\[
 \C[x_1, \dots, x_n]_{S_n} = \C[x_1, \dots, x_n]/
 \langle \C[x_1, \dots, x_n]^{S_n}_{+} \rangle_{\C[x_1, \dots, x_n]}.
\]
We call this algebra the symmetric coinvariant algebra according to
\cite{Hiller}. It is a classical result that we have
\[
 \C[x_1, \dots, x_n]_{S_n} \simeq \C[S_n]
\]
as an $S_n$-module. 

Next, let $A = \C[x,y]/(xy)$ and let 
$M = \{(x,y) \in \C^2\;|\;xy=0\}$ be the corresponding affine variety.
For any $n \in \N$, $S_n$ acts on $A^{\otimes n}$. Set
\[
 (A^{\otimes n})^{S_n}_{+} = \{P \in (A^{\otimes n})^{S_n} \;\vert\; P(0, \dots, 0) = 0\}
\]
and let $J_n$ be the ideal of $A^{\otimes n}$ 
generated by $(A^{\otimes n})^{S_n}_{+}$, i.e.
\[
 J_n = \langle (A^{\otimes n})^{S_n}_{+} \rangle_{A^{\otimes n}}.
\]

\begin{definition}
 The symmetric coinvariant algebra $R_n = A^{\otimes n}_{S_n}$ is
 \[
  R_n = A^{\otimes n}_{S_n} = A^{\otimes n} / J_n.
 \]
\end{definition}

Let $\pi$ be the projection 
\[
\pi : A^{\otimes n} \longrightarrow R_n.
\]
In this paper, we study the
$S_n$-module structure of $R_n$.

By a theorem of Weyl \cite{Weyl}, the elements 
$\sum_{i=1}^{n} x_i^r y_i^s$ 
($r$, $s \geq 0$) generate $\C[x_1, \dots, x_n, y_1, \dots, y_n]^{S_n}$.
In $A^{\otimes n}$, $\sum_{i=1}^{n} x_i^r y_i^s = 0$ for $(r,s)$ such that 
$r \geq 1$ and $s \geq 1$.
Therefore the ideal 
$J_n$ is generated by the power sums $\sum_{i=1}^n x_i^r$ and 
$\sum_{i=1}^n y_i^r$ ($r \geq 1$), or the elementary symmetric
polynomials $e_i$, $f_i$ ($1 \leq i \leq n$).

\section{Generalization of symmetric coinvariant algebra}
\label{sec:some-gener-symm}
In this section, we introduce a generalization of 
$R_n$ and determine its $S_n$-structure.

\begin{definition}
 \label{def:1}
 For $1 \leq i$, $j \leq n$,
let $R^n_{i,j}$ be the quotient algebra of $A^{\otimes n}$ given by
 \[
  R^n_{i,j} = A^{\otimes n} / I^n_{i,j},
 \]
 where
 \begin{eqnarray*}
  S^n_{i,j} & = &
   \left\{
    \begin{array}{l}
     e_1, \dots, e_{i-1}, x_I \quad (|I| = i)\\
     f_1, \dots, f_{j-1}, y_J \quad (|J| = j)\\
    \end{array}
   \right\} \subset A^{\otimes n}, \\
  I^n_{i,j} & = & \langle S^n_{i,j} \rangle_{A^{\otimes n}}.
 \end{eqnarray*}
\end{definition}

Let $\pi^n_{i,j}$ be the
projection
\[
 \pi^n_{i,j} : A^{\otimes n} \longrightarrow R^n_{i,j}.
\]
Clearly, $R_n$ is equal to $R^n_{n,n}$. First we show the following
variant of Newton identity for the elementary symmetric polynomials.

\begin{lemma}[nonsymmetric Newton identity]
 \label{lemma:3}
 For $1 \leq i \leq n$, we have the following identity:
\begin{equation}
 \label{eq:9}
  x_n^i - x_n^{i-1} e_1^{(n)} + \dots + (-1)^{i-1} x_n e_{i-1}^{(n)}
 + (-1)^i e_{i}^{(n)} = (-1)^i e_i^{(n-1)}.
\end{equation}
 \begin{proof}
  Clearly, we have
  \begin{equation}
   \label{eq:8}
    (1+x_n t)^{-1} \prod_{j=1}^{n} (1 + x_j t) = 
    \prod_{j=1}^{n-1} (1 + x_j t).
  \end{equation}
  For $1 \leq i \leq n$, the coefficient of $t^i$ in \refeq{eq:8} 
  coincides with
  \refeq{eq:9} (up to sign).
 \end{proof}
\end{lemma}
The following lemma is easy to prove.
\begin{lemma}
 \label{lemma:5}
 We have equalities $x_n^i = 0$ and $y_n^j = 0$ in $R^n_{i,j}$.
 \begin{proof}
  By \reflemma{lemma:3}, we have
  \[
   x_n^i - x_n^{i-1} e_1^{(n)} + \dots + (-1)^{i-1} x_n e_{i-1}^{(n)}
  + (-1)^{i}
  e_{i}^{(n)} = (-1)^i e_i^{(n-1)} .
  \]
  Since the elements $e_1^{(n)}$, $\dots$, $e_i^{(n)}$, and 
  $e_i^{(n-1)}$ belong
  to $I^n_{i,j}$, we have $x_n^i \in I^n_{i,j}$.
  Similarly, we have $y_n^j \in I^n_{i,j}$.
 \end{proof}
\end{lemma}

In the rest of this section, we determine the $S_n$-module
structure
of $R^n_{i,j}$ for $i+j \leq n+1$. Our proof is a modification of
that in $\cite{GP}$.
First, we introduce another $S_n$-module $R_W$.

For $i$, $j \geq 1$, $i+j \leq n+1$,
let $a_1$, $\dots$, $a_{i-1} \in \C^{\times}$ be distinct, and
let $b_1$, $\dots$, $b_{j-1} \in \C^{\times}$ be also distinct.
We set 
\[
 z_0 = 
 \biggl(
 \binom{a_1}{0}, \dots, \binom{a_{i-1}}{0}, 
 \binom{0}{b_1}, \dots, \binom{0}{b_{j-1}}, 
 \underbrace{\binom{0}{0}, \dots,  \binom{0}{0}}_{n-i-j+2}
 \biggr) \in M^n .
\]
The symmetric group $S_n$ acts on $M^n$.
Let $W$ be the $S_n$-orbit of $z_0$, then $\# W 
= \#(S_n/S_{n-i-j+2})= n!/(n-i-j+2)!$.

We define 
\begin{eqnarray*}
 I_W & = & 
  \left\{ P \in A^{\otimes n} | P(z) = 0 \quad (z \in W) \right\}, \\
 R_W & = & A^{\otimes n} / I_W.
\end{eqnarray*}
The algebra $R_W$ is the coordinate ring of $W \simeq S_n/S_{n-i-j+2}$. Hence 
$R_W  \simeq \Ind^{S_n}_{S_{n-i-j+2}} L_{(n-i-j+2)}$ as an $S_n$-module.

The algebra $A^{\otimes n}$ is graded with the homogeneous degree in
$x$ and $y$. We define the increasing filtration 
$\{G_p A^{\otimes n}\}_{p\geq 0}$ of $A^{\otimes n}$:
$G_p A^{\otimes n}$ is the set of elements of $A^{\otimes n}$ whose
 homogeneous degree are less than $p$.
We also define the filtration $\{G_p R_W\}_{p \geq 0}$ of quotient
algebra $R_W$ as its induced filtration.

\begin{lemma}
 \label{lemma:1}
 Each element of $S^n_{i,j}$ is the leading homogeneous
 component of a polynomial in $I_W$
 \begin{proof}
  For $e_k$ $(1\leq k \leq i-1)$ or $f_l$ $(1 \leq l \leq j-1)$,
  the elements
  \[
   e_k - e_k(a_1, \dots, a_{i-1}, 0, \dots, 0), \quad
  f_l - f_l(b_1, \dots, b_{j-1}, 0, \dots, 0)
  \]
  belong to $I_W$ and their leading homogeneous components are
  $e_k$ or $f_l$.

  The remaining generators $x_I$ $(|I|=i)$ and $y_J$ $(|J|=j)$
  clearly belong to $I_W$. 
 \end{proof}
\end{lemma}

From this lemma, we get the following surjective
homomorphism of $S_n$-modules:
\[
 R^n_{i,j} \longrightarrow \gr R_W
\]
where $\gr R_W$ is the graded algebra associated with the filtration
$\{G_p R_W\}_{p \geq 0}$.
Since the filtration of $R_W$ is invariant by the action of $S_n$,
\reflemma{lemma:6} implies that 
$\gr R_W$ is isomorphic to $R_W$ as an $S_n$-module. Thus we obtain
the following proposition.

\begin{proposition}
 \label{prop:1}
 For $i$, $j \geq 1$ such that $i+j \leq n+1$, 
 there is a surjective homomorphism of $S_n$-modules$:$
 \[
  R^n_{i,j} \longrightarrow \Ind^{S_n}_{S_{n-i-j+2}} L_{(n-i-j+2)}.
 \]
\end{proposition}

Note that $\Dim R^n_{i,j} \geq n!/(n-i-j+2)!$ by \refprop{prop:1}.
Next, we show $\Dim R^n_{i,j} \leq n!/(n-i-j+2)!$. First, we
consider the case of $i+j = n+1$.

We introduce the following filtration of $R^n_{i,j}$ for $i$, 
$j\geq 1$ such that $i+j=n+1$:
\begin{multline*}
 0 = \langle y_n^j \rangle_{R^n_{i,j}} \subset 
 \langle y_n^{j-1} \rangle_{R^n_{i,j}} \subset \dots \subset
 \langle y_n^{1} \rangle_{R^n_{i,j}} \\
 \subset 
 \langle y_n, \X^{(n-1)}_{i-1} \rangle_{R^n_{i,j}}
 \subset 
 \langle x_n^{i-2}, y_n, \X^{(n-1)}_{i-1} \rangle_{R^n_{i,j}} 
 \dots 
 \subset
 \langle x_n^1, y_n, \X^{(n-1)}_{i-1} \rangle_{R^n_{i,j}} 
 \subset R^n_{i,j}
\end{multline*}
where $\X^{(n-1)}_{i-1} = \{ x_I \;|\; I \subset \{1, \dots, n-1\}, |I|=i-1\}$.
Note that, since $x_n^{i-1} \equiv (-1)^{i-1} e^{(n-1)}_{i-1}$ 
by \reflemma{lemma:3},
 $\langle y_n, \X^{(n-1)}_{i-1} \rangle_{R^n_{i,j}} =
\langle y_n, x_n^{i-1}, \X^{(n-1)}_{i-1} \rangle_{R^n_{i,j}}$.
From this filtration, we have the following decomposition of 
$R^n_{i,j}$:
\begin{multline}
 \label{eq:4}
 R^n_{i,j} \simeq R^n_{i,j} / 
 \langle x_n, y_n, \X^{(n-1)}_{i-1} \rangle_{R^n_{i,j}} \\
 \oplus 
 \bigoplus_{k=1}^{i-1} \langle x_n^k, y_n, \X^{(n-1)}_{i-1} \rangle_{R^n_{i,j}} /
 \langle x_n^{k+1}, y_n, \X^{(n-1)}_{i-1} \rangle_{R^n_{i,j}} \\
 \oplus
 \bigoplus_{k=1}^{j-1} \langle y_n^k \rangle_{R^n_{i,j}} /
 \langle y_n^{k+1} \rangle_{R^n_{i,j}}.
\end{multline}

\begin{lemma}
 \label{lemma:8}
 For $i$, $j \geq 1$ such that $i+j = n+1$, $1 \leq k \leq i-1$ and 
$1 \leq k' \leq j-1$, we have the following surjective 
homomorphisms of $S_{n-1}$-modules$:$
\begin{align*}
 \varphi_0 : R^{n-1}_{i-1, j} &\rightarrow 
 R^n_{i,j} / \langle y_n, x_n, \X^{(n-1)}_{i-1} \rangle_{R^n_{i,j}}, \\
 P + I^{n-1}_{i-1,j} &\mapsto P + I^n_{i,j} + 
 \langle y_n, x_n, \X^{(n-1)}_{i-1} \rangle_{A^{\otimes n}}, \\
 \varphi_k : R^{n-1}_{i-1,j}  &\rightarrow
 \langle y_n, x_n^k, \X^{(n-1)}_{i-1} \rangle_{R^n_{i,j}}/
 \langle y_n, x_n^{k+1}, \X^{(n-1)}_{i-1} \rangle_{R^n_{i,j}}, \\
 P + I^{n-1}_{i-1,j} &\mapsto x_n^k P + I^n_{i,j} + 
 \langle y_n, x_n^{k+1}, \X^{(n-1)}_{i-1} \rangle_{A^{\otimes n}}, \\
 \varphi'_{k'} : R^{n-1}_{i,j-1}&\rightarrow
 \langle y_n^k \rangle_{R^n_{i,j}}/
 \langle y_n^{k+1} \rangle_{R^n_{i,j}}, \\
 P + I^{n-1}_{i,j-1} &\mapsto y_n^k P + I^n_{i,j} + 
 \langle y_n^{k+1} \rangle_{A^{\otimes n}}.
\end{align*}
 \begin{proof}
  First, since $e^{(n-1)}_l \equiv e^{(n)}_l$ modulo 
$\langle x_n \rangle_{A^{\otimes n}}$ and $f^{(n-1)}_l \equiv f^{(n)}_l$
modulo $\langle y_n \rangle_{A^{\otimes n}}$, we have 
  $I^{n-1}_{i-1,j} \subset I^n_{i,j} + \langle y_n, x_n,
  \X^{(n-1)}_{i-1}\rangle_{A^{\otimes n}}$. 
  Therefore, $\varphi_0$ is well-defined, and
  clearly it is surjective. 

  Next, we show that $\varphi_k$ is well-defined for $1 \leq k \leq
  i-1$. Let $P \in A^{\otimes n-1} \subset A^{\otimes n}$, and assume
  $P$ belongs to $I^{n-1}_{i-1,j}$. We have
  \begin{multline*}
   P = P_1 e^{(n-1)}_1 + \dots + P_{i-2} e^{(n-1)}_{i-2} +
   \sum_{\substack{I \subset \{1, \dots, n-1\} \\ |I|=i-1}} P_I x_I \\
   + Q_1 f^{(n-1)}_1 + \dots + Q_{j-1} f^{(n-1)}_{j-1}
   + \sum_{\substack{J \subset \{1, \dots, n-1\} \\ |J|=j}} Q_J y_J
  \end{multline*}
  where $P_1$, $\dots$, $P_{i-2}$, $P_I$, $Q_1$, $\dots$, $Q_{j-1}$,
  $Q_J \in A^{\otimes n-1}$. Therefore, we have
  \begin{multline*}
   x_n^k P = P_1 x_n^k e^{(n-1)}_1 + \dots + P_{i-2} x_n^k e^{(n-1)}_{i-2} +
   \sum_{\substack{I \subset \{1, \dots, n-1\} \\ |I|=i-1}} x_n^k P_I x_I \\
   +Q_1 x_n^k f^{(n-1)}_1 + \dots + Q_{j-1} x_n^k f^{(n-1)}_{j-1}
   +\sum_{\substack{J \subset \{1, \dots, n-1\} \\ |J|=j}} x_n^k Q_J y_J.
  \end{multline*}
  Since $x_n^k e^{(n-1)}_l \equiv x_n^k e^{(n)}_l$ modulo 
  $\langle x_n^{k+1} \rangle_{A^{\otimes n}}$ and
  $x_n^k f^{(n-1)}_{l} = x_n^k f^{(n)}_l$, $x_n^k P$ belongs to
  $I^n_{i,j} + \langle y_n, x_n^{k+1}, \X^{(n-1)}_{i-1} 
  \rangle_{A^{\otimes n}}$. Therefore,
  $\varphi_k$ is well-defined.

  For any element of 
  $\langle y_n, x_n^k, \X^{(n-1)}_{i-1} \rangle_{R^n_{i,j}} /
  \langle y_n, x_n^{k+1}, \X^{(n-1)}_{i-1} \rangle_{R^n_{i,j}}$, we
  can choose its representative $x_n^k P$ where $P \in A^{\otimes n-1}$.
  Therefore $\varphi_k$ is surjective.

  Similarly, $\varphi'_{k'}$ is well-defined and surjective.
 \end{proof}
\end{lemma}

\begin{proposition}
 \label{prop:2}
 For $i$, $j \geq 1$ such that $i+j = n+1$, the $S_n$-module
 $R^n_{i,j}$ is isomorphic to the regular representation of
 $S_n$.
 \begin{proof}
  By \refprop{prop:1}, we have the surjective homomorphism
  $R^n_{i,j} \twoheadrightarrow \C[S_n]$. We show
  $\Dim R^n_{i,j} \leq n!$ by induction on $n$. First, 
  consider the case of $n=1$. In this case, we have $i=j=1$,
  and this case is already proved. 
  We may assume that $\Dim R^{n-1}_{i',j'} = (n-1)!$ for $i'$, 
  $j' \geq 1$ such that $i'+j'=n$. Therefore,
  by \refeq{eq:4} and
  \reflemma{lemma:8}, we have
  \begin{align*}
   \Dim R^n_{i,j} &\leq \Dim R^{n-1}_{i-1,j} +
   \sum_{k=1}^{i-1} \Dim R^{n-1}_{i-1,j} +
   \sum_{k'=1}^{j-1} \Dim R^{n-1}_{i,j-1} \\
   &= n!.
  \end{align*}
  Hence, the induction completes.
 \end{proof}
\end{proposition}

Next, consider the case of $i+j \leq n$. We introduce the following
filtration of $R^n_{i,j}$ for $i$, $j \geq 1$ such that
$i+j \leq n$:
\begin{multline*}
  0 = \langle y_n^j \rangle_{R^n_{i,j}} \subset 
  \langle y_n^{j-1} \rangle_{R^n_{i,j}} \subset \dots \subset
 \langle y_n^{1} \rangle_{R^n_{i,j}} \\
 \subset 
 \langle x_n^{i-1}, y_n \rangle_{R^n_{i,j}}
 \subset 
 \langle x_n^{i-2}, y_n \rangle_{R^n_{i,j}} \subset
 \dots 
 \subset
 \langle x_n^1, y_n \rangle_{R^n_{i,j}}
 \subset
 R^n_{i,j}.
\end{multline*}
From this filtration, we have the following decomposition of 
$R^n_{i,j}$:
\begin{multline}
 \label{eq:10}
 R^n_{i,j} \simeq R^n_{i,j} / 
 \langle x_n, y_n \rangle_{R^n_{i,j}}
 \oplus 
 \bigoplus_{k=1}^{i-1} \langle x_n^k, y_n \rangle_{R^n_{i,j}} /
 \langle x_n^{k+1}, y_n \rangle_{R^n_{i,j}} \\
 \oplus
 \bigoplus_{k=1}^{j-1} \langle y_n^k \rangle_{R^n_{i,j}} /
 \langle y_n^{k+1} \rangle_{R^n_{i,j}}.
\end{multline}

We can prove the following lemma similarly to \reflemma{lemma:8}.
\begin{lemma}
 \label{lemma:9}
 For $i$, $j \geq 1$ such that $i+j \leq n$, $1 \leq k \leq i-1$ and 
 $1 \leq k' \leq j-1$, we have the following surjective 
homomorphisms of $S_{n-1}$-modules$:$
\begin{align*}
 \varphi_0 : R^{n-1}_{i, j} &\rightarrow 
 R^n_{i,j} / \langle y_n, x_n \rangle_{R^n_{i,j}}, \\
 P + I^{n-1}_{i,j} &\mapsto P + I^n_{i,j} + 
 \langle y_n, x_n \rangle_{A^{\otimes n}}, \\
 \varphi_k : R^{n-1}_{i-1,j}  &\rightarrow
 \langle y_n, x_n^k \rangle_{R^n_{i,j}}/
 \langle y_n, x_n^{k+1} \rangle_{R^n_{i,j}}, \\
 P + I^{n-1}_{i-1,j} &\mapsto x_n^k P + I^n_{i,j} + 
 \langle y_n, x_n^{k+1} \rangle_{A^{\otimes n}}, \\
 \varphi'_{k'} : R^{n-1}_{i,j-1}&\rightarrow
 \langle y_n^k \rangle_{R^n_{i,j}}/
 \langle y_n^{k+1} \rangle_{R^n_{i,j}}, \\
 P + I^{n-1}_{i,j-1} &\mapsto y_n^k P + I^n_{i,j} + 
 \langle y_n^{k+1} \rangle_{A^{\otimes n}}.
\end{align*}
\end{lemma}

\begin{proposition}
 \label{prop:3}
 For $i$, $j \geq 1$ such that $i+j \leq n$, we have the
 following $S_n$-module isomorphism$:$
 \[
 R^n_{i,j} \simeq \Ind^{S_n}_{S_{n-i-j+2}} L_{(n-i-j+2)}.
 \]
 \begin{proof}
  The following proof is similar to one of \refprop{prop:2}.
  By \refprop{prop:1}, we have the surjective homomorphism
  $R^n_{i,j} \twoheadrightarrow \Ind^{S_n}_{S_{n-i-j+2}} L_{(n-i-j+2)}$.
  We show
  $\Dim R^n_{i,j} \leq n!/(n-i-j+2)!$ by induction on $n$. First, 
  consider the case of $n=2$. In this case, we have $i=j=1$,
  and this case is already proved. We may
  assume that $\Dim R^{n-1}_{i',j'} = (n-1)!/(n-i'-j'+1)!$ for $i'$, 
  $j' \geq 1$ such that $i'+j'\leq n-1$. By \refprop{prop:2},
  we have that
  $\Dim R^n_{i,j} = n! = n!/(n-i-j+2)!$ for $i$, $j \geq 1$ such that 
  $i+j = n+1$. Therefore,
  by \refeq{eq:10} and
  \reflemma{lemma:9}, we have
  \begin{align*}
   \Dim R^n_{i,j} &\leq \Dim R^{n-1}_{i,j} +
   \sum_{k=1}^{i-1} \Dim R^{n-1}_{i-1,j} +
   \sum_{k'=1}^{j-1} \Dim R^{n-1}_{i,j-1} \\
   &= \frac{(n-1)!}{(n-i-j+1)!} + (i-1) \frac{(n-1)!}{(n-i-j+2)!}
   + (j-1) \frac{(n-1)!}{(n-i-j+2)!} \\
   &= \frac{n!}{(n-i-j+2)!}.
  \end{align*}
  Therefore, the induction completes.
 \end{proof}
\end{proposition}

\section{The structure of $R_n$}
\label{sec:structure-r_n}
In this section, we determine the $S_n$-module structure of $R_n$.

We define a decreasing filtration 
$\{F^i A^{\otimes n}\}_{0 \leq i \leq n}$ of $A^{\otimes n}$ where
\[
 F^i A^{\otimes n} = \sum_{|J|=i,\:J \subset \{1, \dots, n\}} y_J A^{\otimes n}.
\]
This filtration is $S_n$-invariant. 
Let $F^i J_n = J_n \cap F^i A^{\otimes n}$ and
$F^i R_n = \pi(F^i A^{\otimes n})$.

Let $R_n^{(i)} = gr^i R_n = F^i R_n / F^{i+1} R_n = F^i A^{\otimes n} / 
\left(F^i J_n + F^{i+1} A^{\otimes n}\right)$, then $R_n^{(n)} = 0$.
We have
\begin{equation}
 \label{eq:18}
 R_n^{(0)} = A^{\otimes n} / \left(F^0 J_n + F^1 A^{\otimes n}\right) =
 \C[x_1, \dots, x_n]_{S_n} \simeq \C[S_n]
\end{equation}
as an $S_n$-module (\cite{Chev}).

Since the algebra $R_n$ has the $S_n$-invariant filtration $\{F^i
R_n\}_{0\leq i \leq n}$, 
$R_n$ is isomorphic to $\gr R_n = \bigoplus_{i=0}^{n-1} R_n^{(i)}$
by \reflemma{lemma:6}. 
For $1 \leq i \leq n-1$, we will determine the $S_n$-module
 structure of 
$R_n^{(i)}$ by using the result of \refsec{sec:some-gener-symm}.
Since $F^i J_n + F^{i+1} A^{\otimes n} \subset I^n_{n-i,i+1}$, there is
a homomorphism $\phi: R_n^{(i)} \rightarrow R^n_{n-i,i+1}$.

\begin{definition}
 Let $\A$ be a commutative ring and let $\M$ be an $\A$-module.
 \begin{enumerate}
  \item An element $a \in \A$ is called $\M$-regular if and only if 
	for any $0 \ne x \in \M$ we have $a x \ne 0$.
  \item A sequence $a_1$, $\dots$, $a_n \in \A$ is called 
	an $\M$-regular sequence
	if and only if for $j = 1$, $\dots$, $n$, $a_j$ is 
	$(\M / \sum_{k=1}^{j-1} a_k \M)$-regular.
 \end{enumerate}
\end{definition}

\begin{lemma}
 \label{lemma:4}
 For any $n \in \N$, the sequence of the elementary symmetric polynomials
 $e_1$, $\dots$, $e_n$ is the $\C[x_1, \dots, x_n]$-regular sequence.
\end{lemma}

\begin{lemma}
 \label{lemma:2}
 Let $\M$ be a flat $\A$-module.
 If $f_1$, $\dots$, $f_n \in \A$ is an $\A$-regular sequence,
 $f_1$, $\dots$, $f_n$ is an $\M$-regular
 sequence.
\end{lemma}

These lemmas are basic facts in the theory of commutative algebra.

\begin{proposition}
 \label{prop:4}
 For $1 \leq i \leq n-1$,
 the homomorphism $\phi: R_n^{(i)} \rightarrow R^n_{n-i,i+1}$ is 
 injective.
 \begin{proof}
  By the definition of $\phi$, we only need to prove that 
  $J_n + F^{i+1} A^{\otimes n} \supset F^i I^n_{n-i,i+1}$ for $1 \leq i \leq n-1$.

  For $J \subset \{1, \dots, n\}$, let $\bar{J}$ be the complement 
  of $J$ in $\{1, \dots, n\}$.
  Fix an arbitrary element $P \in F^i I^n_{n-i,i+1}$. We can decompose
  $P$ into two forms. First,
  \begin{equation}
   \label{eq:14}
   P = \sum_{|J|=i} P_J y_J + P'
  \end{equation}
  where $P_J \in \C[x_j (j \in \bar{J})] \otimes \C[y_j (j \in J)]$
  and $P' \in F^{i+1} A^{\otimes n}$. Second,
  \begin{equation}
   \label{eq:15}
   P = Q_1 e_1 + \dots + Q_{n-i-1} e_{n-i-1} + \sum_{|I|=n-i} Q_I x_I
    + R_1 f_1 + \dots + R_{i} f_i + \sum_{|J|=i+1} R_J y_J
  \end{equation}
  where $Q_1$, $\dots$, $Q_{n-i-1}$, $Q_I$ $R_1$, $\dots$, $R_{i}$,
  $R_J \in A^{\otimes n}$. Fix $J \subset \{1, \dots, n\}$, $|J|=i$.
  For $S \in A^{\otimes n}$,
  we denote by $\bar{S}$ the element of
  $\C[x_j (j \in \bar{J})] \otimes \C[y_j (j \in J)]$
  obtained from the substitution $x_j = 0$ ($j\in J$) and
  $y_j = 0$ ($j \in \bar{J}$).
  Set $x_j = 0$ for $j \in J$ and $y_j = 0$ for 
  $j \in \bar{J}$ in \refeq{eq:14} and \refeq{eq:15},
  we have
  \begin{multline*}
   P_J y_J =
   \bar{Q}_1 e_1(\{x_j\}_{j\in\bar{J}}) + \dots +
   \bar{Q}_{n-i-1} e_{n-i-1}(\{x_j\}_{j\in\bar{J}}) + 
   \bar{Q}_{\bar{J}} x_{\bar{J}} \\
   + \bar{R}_1 f_1(\{y_j\}_{j\in J}) + \dots +
   + \bar{R}_{i-1} f_{i-1}(\{y_j\}_{j\in J}) + 
   \bar{R}_{i} y_{J}.
  \end{multline*}
  Since $x_{\bar{J}} = e_{n-i}(\{x_j\}_{j \in \bar{J}})$ and
  $y_J = f_i(\{y_j\}_{j\in J})$, we have
  \begin{equation}
 \label{eq:5}
 (P_J - \bar{R}_i) f_i(\{y_j\}_{j\in J})
  \in \left\langle 
   \!\!
  \begin{array}{l}
  e_1(\{x_j\}_{j\in\bar{J}}), \dots,
  e_{n-i}(\{x_j\}_{j\in\bar{J}}), 
   \\
  f_1(\{y_j\}_{j\in J}), \dots, f_{i-1}(\{y_j\}_{j\in J})
  \end{array}
  \!\!
  \right\rangle_{\C[x_j (j \in \bar{J})] \otimes \C[y_j (j \in J)]}
  \end{equation}
  in $\C[x_j (j \in \bar{J})] \otimes \C[y_j (j \in J)]$. Let
  $\A = \C[y_j (j \in J)]$ and let 
  $\M = \C[x_j (j \in \bar{J})]_{S_{n-i}} \otimes \C[y_j (j \in J)]$.
  By \reflemma{lemma:4},
  $f_1(\{y_j\}_{j\in J})$, $\dots$, $f_{i}(\{y_j\}_{j\in J})$ is
  the $\A$-regular sequence, and $\M$
  is a flat $\A$-module. 
  Hence, by \reflemma{lemma:2},
  $f_1(\{y_j\}_{j\in J})$, $\dots$, $f_{i}(\{y_j\}_{j\in J})$ is
  an $\M$-regular sequence. Therefore, from \refeq{eq:5}, we have
  \[
   P_J - \bar{R}_i
  \in \left\langle 
  \begin{array}{l}
  e_1(\{x_j\}_{j\in\bar{J}}), \dots,
  e_{n-i}(\{x_j\}_{j\in\bar{J}}), 
   \\
  f_1(\{y_j\}_{j\in J}), \dots, f_{i-1}(\{y_j\}_{j\in J})
  \end{array}
  \right\rangle_{\C[x_j (j \in \bar{J})] \otimes \C[y_j (j \in J)]}.
  \]
  By multiplying $y_J$ we have
  \[
  (P_J - R_i) y_J
  \in \left\langle 
  \begin{array}{l}
  e_1, \dots, e_{n-i}, 
   \\
  f_1, \dots, f_{i-1}
  \end{array}
  \right\rangle_{A^{\otimes n}} + F^{i+1} A^{\otimes n}.
  \]
  Therefore, taking the summation for $J$ ($|J|=i$), we have
  \[
   P - P' - R_i f_i 
  \in \left\langle 
  \begin{array}{l}
  e_1, \dots, e_{n-i}, 
   \\
  f_1, \dots, f_{i-1}
  \end{array}
  \right\rangle_{A^{\otimes n}} + F^{i+1} A^{\otimes n} \subset J_n + F^{i+1} A^{\otimes n}.
  \]
  Hence $P \in J_n + F^{i+1} A^{\otimes n}$.
 \end{proof}
\end{proposition}

We have shown that $\phi: R_n^{(i)} \rightarrow R^n_{n-i,i+1}$ is
an injective homomorphism.
On the other hand, we have a surjective homomorphism 
$\psi: R^n_{n-i,i+1} \rightarrow R^n_{n-i,i}$ because 
$I^n_{n-i,i+1} \subset I^n_{n-i,i}$.

\begin{proposition}
 \label{prop:5}
 For $1 \leq i \leq n-1$, we have the following exact sequence of
 $S_n$-modules.
 \[
 0 \rightarrow R_n^{(i)} \xrightarrow{\phi} R^n_{n-i,i+1}
 \xrightarrow{\psi} R^n_{n-i,i} \rightarrow 0
 \]
 \begin{proof}
  First we have $\Im \phi \subset \Ker \psi$ because 
  $F^{i} A^{\otimes n} \subset I^n_{n-i,i}$. 

  Next we show that $\Im \phi \supset \Ker \psi$. If
  $P \in \Ker \psi$, it belongs to the image of $I^n_{n-i,i}$
  in $R^n_{n-i,i+1}$. Namely we have  in $R^n_{n-i,i+1}$
  \begin{eqnarray*}
   P & = & P_1 e_1 + \dots + P_{n-i-1} e_{n-i-1} + 
    \sum_{|I|=n-i} P_I x_I \\
   & & \hspace*{1.5cm}+ Q_1 f_1 + \dots + Q_{i-1} f_{i-1} +
    \sum_{|J|=i} Q_J y_J \\
   & = & \sum_{|J|=i} Q_J y_J \in F^i R^n_{n-i,i+1} = \phi(R_n^{(i)}).
  \end{eqnarray*}
 \end{proof}
\end{proposition}

\begin{corollary}
 \label{cor:3}
 For $1 \leq i \leq n-1$, we have the following isomorphism
 of $S_n$-modules$:$
 \[
  R_n^{(i)} \simeq \Ind^{S_n}_{S_2} L_{(1,1)}.
 \]
 \begin{proof}
  By \refprop{prop:2} and \refprop{prop:3}, 
  $R^n_{n-i,i+1} \simeq \C[S_n]$ and
  $R^n_{n-i,i} \simeq \Ind^{S_n}_{S_2} L_{(2)}$, the claim
  of the corollary follows from the exact sequence of \refprop{prop:5}.
 \end{proof}
\end{corollary}

Together with \refeq{eq:18}, we obtain the following theorem.
\begin{theorem}
 \label{thm:1}
 We have the following isomorphism of $S_n$-modules$:$
 \[
  R_n \simeq \C[S_n] \oplus (n-1) \Ind^{S_n}_{S_2} L_{(1,1)}.
 \]
\end{theorem}

\section{The local Weyl module at a double point}
\label{sec:local-weyl-module}
In this section, we study the structure
of the local Weyl module at the double point.

Let $\g$ be a finite-dimensional semisimple Lie algebra and let
$\g = \mathfrak{n}_{+} \oplus \mathfrak{h} \oplus \mathfrak{n}_{-}$ be
its triangular decomposition. Let $M$ be an affine variety and let $A$
be the coordinate ring of $M$.
In \cite{FL}, Feigin and Loktev introduced the $(\g \otimes A)$-module
$W_M(\{0\}_{\lambda})$ called the local Weyl module for a dominant
integrable weight $\lambda \in \mathfrak{h}^{*}$.
$W_M(\{0\}_{\lambda})$ is the maximal
$\g$-integrable module with a cyclic vector $v_0$ such that:
\[
 (\mathfrak{n}_{+} \otimes P) v_0 = 0, \qquad
 (h \otimes P) v_0 = \lambda(h) P(0) v_0 \quad (P \in A, h \in \mathfrak{h}).
\]

Consider the case of  $\g = \Slr$. Let
$V_{r+1}$ be the vector representation of $\Slr$ and let
$\omega_1$ be the highest weight of $V_{r+1}$.
For $\lambda = n \omega_1$,
the following theorem is proved by Feigin and Loktev.
\begin{theorem}[\cite{FL}]
 \label{thm:2}
 There is an isomorphism of $\Slr$-modules$:$
 \[
  W_M(\{0\}_{n \omega_1}) \simeq 
 \left(V_{r+1}^{\otimes n} \otimes A^{\otimes n}_{S_n} 
 \right)^{S_n}.
 \]
\end{theorem}

Thus, combining \refthm{thm:1} and \refthm{thm:2},
we obtain the $\Slr$-module structure of $W_M(\{0\}_{n \omega_1})$ as follows.

\begin{proposition}
 \label{prop:6}
 For $n \in \Z_{\geq 0}$, we have the following isomorphism of
 $\Slr$-modules.
 \[
  W_M(\{0\}_{n \omega_1}) \simeq V_{r+1}^{\otimes n} \oplus
 (n-1) \left(
 V_{r+1}^{\otimes n-2} \otimes \wedge^{2} V_{r+1}
 \right).
 \]
 \begin{proof}
  The following proof is essentially same as the first half of the
  proof of Theorem~10 in \cite{FL}.

  By \refthm{thm:2} and \refthm{thm:1}, we have
  \begin{align*}
   W_M(\{0\}_{n \omega_1}) 
   &\simeq
   \left(
   V^{\otimes n}_{r+1} \otimes A^{\otimes n}_{S_n}
   \right)^{S_n} \\
   &\simeq
   \left(
   V^{\otimes n}_{r+1} \otimes (\C[S_n] \oplus (n-1) \Ind^{S_n}_{S_2} L_{(1,1)})
   \right)^{S_n} \\
   &\simeq
   \left(
   V^{\otimes n}_{r+1} \otimes \C[S_n]
   \right)^{S_n}
   \oplus (n-1) \left(V^{\otimes n}_{r+1} \otimes
   \Ind^{S_n}_{S_2} L_{(1,1)}
   \right)^{S_n} \\
   &\simeq
   \left(
   V^{\otimes n}_{r+1} \otimes \C[S_n]
   \right)^{S_n}
   \oplus (n-1) \left(V^{\otimes n}_{r+1} \otimes
    L_{(1^n)} \otimes
   \Ind^{S_n}_{S_2} L_{(2)}
   \right)^{S_n} \\
   &\simeq
   V^{\otimes n}_{r+1}
   \oplus (n-1)
   \Hom_{S_n}(L_{(1^n)} \otimes (V_{r+1}^{*})^{\otimes n},
   \Ind^{S_n}_{S_2} L_{(2)}) \\
   &\simeq
   V_{r+1}^{\otimes n}
   \oplus (n-1) 
   \Hom_{S_2}
   (L_{(1^n)} \otimes (V_{r+1}^{*})^{\otimes n},
   L_{(2)}) \\
   &\simeq
   V_{r+1}^{\otimes n} \oplus (n-1)
   V_{r+1}^{\otimes n-2} \otimes \wedge^2 V_{r+1}.
  \end{align*}
 \end{proof}
\end{proposition}

\begin{corollary}
 For $n \in \Z_{\geq 0}$, we have
 \[
  \Dim W_M(\{0\}_{n \omega_1}) = (r+1)^{n-2}
 \left((r+1)^2+\frac{(n-1)(r+1)r}{2} \right).
 \]
\end{corollary}

\end{document}